\documentclass{amsart}

\usepackage{amssymb,mathrsfs}
\usepackage{color,umoline}
\usepackage[dvipsnames]{xcolor}
\usepackage{graphicx}
\usepackage{enumitem}
\usepackage[font=small,labelfont=bf]{caption}
\usepackage{tikz, caption, subcaption}
\usetikzlibrary{quotes,angles}
\usepackage{relsize}
\usepackage[mathscr]{eucal}
\usepackage{verbatim}
\usepackage[draft]{todonotes}
\usepackage{hyperref}
\usepackage[dvipsnames]{xcolor}
\usetikzlibrary{patterns}
\tikzstyle{EDR}=[draw=lightgray,line width=0pt,preaction={clip, postaction={pattern=north east lines, pattern color=gray}}]
\tikzstyle{EDR1}=[draw=lightgray,line width=0pt,preaction={clip, postaction={pattern=north west lines, pattern color=gray}}]

\newtheorem{theorem}{Theorem}[section]
\newtheorem{lemma}[theorem]{Lemma}
\newtheorem{proposition}[theorem]{Proposition}

\theoremstyle{definition}
\newtheorem{definition}[theorem]{Definition}

\theoremstyle{remark}
\newtheorem{remark}[theorem]{Remark}

\numberwithin{equation}{section}

\theoremstyle{remark}

\numberwithin{equation}{section}

\definecolor{mygray}{gray}{0.95}

\definecolor{mypink1}{rgb}{1.2,1.1,0.9}

\definecolor{mypink2}{rgb}{1.0,0.95 ,0.9}

\definecolor{mypink3}{rgb}{1.0,0.6,0.7}

\newcommand{\R}{\mathbb{R}}

\begin{document}
	\title[Decoupling]{Decoupling for convex hypersurfaces of finite type in higher dimensions}
	
	\author[C. Gao]{Chuanwei Gao}
	\address{School of Mathematical Sciences, Capital Normal University, Beijing, 100048, China}
	\email{cwgao@pku.edu.cn}
	
	\author[Z. Li]{Zhuoran Li}
	\address{Department of Mathematics, Taizhou University, Taizhou 225300, China}
	\email{lizhuoran18@gscaep.ac.cn}

	\author[T. Zhao]{Tengfei Zhao}%
	\address{School of Mathematics and Physics, University of Science and Technology Beijing, Beijing	100083, China}
	\email{zhao\underline{ }tengfei@ustb.edu.cn}%

	\author[J. Zheng]{Jiqiang Zheng}
	\address{Institute of Applied Physics and Computational Mathematics, Beijing 100088}
	\email{zhengjiqiang@gmail.com}

	\begin{abstract}
In this paper, we establish an $\ell^2$ decoupling inequality for the convex hypersurface
		\[\Big\{(\xi_1,...,\xi_{n-1},\xi_1^m+...+\xi_{n-1}^m): (\xi_1,...,\xi_{n-1}) \in [0,1]^{n-1}\Big\}\]associated with the decomposition adapted to hypersurfaces of finite type,
where $n\geq 2$ and  $m\geq 4$ is an even number.
The key ingredients of the proof include $\ell^2$ decoupling inequalities for convex hypersurfaces of the form
		\[\Big\{(\xi_1,...,\xi_{n-1},\phi_1(\xi_1)+...+\phi_s(\xi_s)+\xi_{s+1}^m+...+\xi_{n-1}^m): (\xi_1,...,\xi_{n-1}) \in [0,1]^{n-1}\Big\},\]
		$0 \leq s \leq n-1$, with $\phi_1,...,\phi_s$ being $m$-nondegenerate.
	\end{abstract}
	
	

	\maketitle	





\section{Introduction and main result}\label{sec:intmai}

\noindent

Decoupling inequality was introduced by Wolff \cite{Wolff00} with a purpose to study the local smoothing estimates for the solution to the wave equation. After a series of work \cite{GaSe09,GaSchSe,LaWo02}, finally,  Bourgain and Demeter \cite{BD15} proved the sharp $\ell^2$ decoupling inequality for compact hypersurfaces with positive definite second fundamental form and followed by a wide range of important applications \cite{BD15, BDG16, GWZ20}.
 It is also an interesting project to generalize the decoupling inequality to hypersurfaces
with some kind of degeneracy.
 For that direction, we refer to
\cite{Xi, Yang21, Kemp, LiYang21preprint, LiZheng21preprint}.
In this paper, we will study the $\ell^2$ decoupling problem for certain convex hypersurfaces of finite type in $\mathbb{R}^n$.
 More precisely, we consider the hypersurface in $\mathbb{R}^n$ given by
\begin{equation*}
  F^{n-1}_{m}(0,n-1):=\Big\{(\xi_1,...,\xi_{n-1},\xi_1^m+...+\xi_{n-1}^m): (\xi_1,...,\xi_{n-1}) \in [0,1]^{n-1}\Big\},\;m\geq 2.
\end{equation*}
For any function $g\in L^1([0,1]^{n-1})$ and each subset $Q\subset [0,1]^{n-1}$, we denote the corresponding Fourier extension operator by
\begin{equation}\label{eq:001}
  E_{Q}g(x):= \int_Q g(\xi_1,...,\xi_{n-1})e\big(x_1\xi_1+...+x_{n-1}\xi_{n-1}+x_n(\xi^m_1+...+\xi_{n-1}^m)\big)\;d\xi_1...d\xi_{n-1},
\end{equation}
where $e(t)=e^{2\pi i t}$ for $t \in \mathbb{R}$, and $x=(x_1,...,x_n)\in\R^n.$

For $m=2$, the hypersurface $F^{n-1}_{2}(0,n-1)$ is exactly a paraboloid over the region $[0,1]^{n-1}$. For $n=3,\;m=4$, the hypersurface $F^2_{ 4}(0,2)$ is exactly the surface studied in  \cite{LiZheng21preprint}.
In terms of $m=2$, Bourgain-Demeter\cite{BD15} showed that
\begin{theorem}\label{thm:BD}
Let $S:=\{\xi,\psi(\xi)\} \subset \R^n$ be a smooth hypersurface with positive definite second fundamental form and $E^S_{[0,1]^{n-1}}$ be an extension operator defined as above
associated with the graph of $\psi$.
Suppose that $2\leq p\leq \frac{2(n+1)}{n-1}$.
 For each $\varepsilon > 0$,
 there exists a constant $C(\varepsilon,p)$
 such that
\begin{equation}\label{equ:bd} \|E_{[0,1]^{n-1}}^Sg\|_{L^{p}(B_R)}\leq C(\varepsilon,p)R^{\varepsilon}\Big(\sum_{\delta\in \mathbf{q}} \|E_{\delta}^S g\|^2_{L^{p}(w_{B_R})}\Big)^{1/2},
\end{equation}
where $w_{B_R}(x)= \big(1+\frac{\vert x-x_0 \vert}{R}\big)^{-100n}$
denotes the standard weight function adapted to the ball $B_R$ in $\mathbb{R}^n$ centered at $x_0$ with radius $R$ and $\mathbf{q}$ is a family of finitely overlapping cubes of dimension $R^{-1/2}$
such that
\begin{equation*}
[0,1]^{n-1}=\bigcup_{\delta\in \mathbf{q}} \delta.
\end{equation*}
\end{theorem}

For convenience, we use $P^{n-1}$ to denote $F_{2}^{n-1}(0,n-1)$.
Suppose that the Fourier support of $F$ is contained in the $\frac{1}{R}$-neighborhood of $P^{n-1}$ which we denote by $\mathcal{N}_{1/R}(P^{n-1})$.
We decompose
\[
\mathcal{N}_{1/R}(P^{n-1})
=
\bigcup_{\delta \in \mathbf{q}}
\bar{\delta},
\]
where $\bar{\delta}$ denotes the $\frac{1}{R}$-neighborhood of $\{(\xi,\psi(\xi)): \xi\in\delta\}$. Roughly speaking, each $\bar{\delta}$ can be seen as a rectangular box of dimension $R^{-1/2}\times \cdots\times R^{-1/2}\times R^{-1}$.
Correspondingly, we decompose $F$ as
\begin{equation*}
F=\sum_{\bar{\delta}} F_{\bar{\delta}},
\end{equation*}
where $F_{\bar{\delta}}:= \mathcal{F}^{-1}\big(\hat{F}\chi_{\bar{\delta}}\big)$ and $\chi_{\bar{\delta}}$ is the characteristic function on $\bar{\delta}$.
An equivalent form of \eqref{equ:bd} is given by
\begin{equation}\label{equ:bd1}
  \|F\|_{L^{p}(\mathbb{R}^n)}\leq C(\varepsilon,p)R^{\varepsilon}\Big(\sum_{\bar{\delta}}\|F_{\bar{\delta}}\|^2_{L^{p}(\mathbb{R}^n)}\Big)^{1/2},\;\; 2\leq p \leq \frac{2(n+1)}{n-1}.
\end{equation}

For $m>2$, the Gaussian curvature of the hypersurface $F^{n-1}_{m}(0,n-1)$ vanishes at the points $\xi$ where there is at least one $j$ such that $\xi_j=0$.
 In the following, we assume that $m\geq 4$ is an even number. Given $R\gg1$, we divide $[0,1]$ into
\begin{equation*}
  [0,1]=\bigcup_k I_{k}.
\end{equation*}
where $I_{0}=[0,R^{-\frac1m}]$ and
$$I_{k}=[2^{k-1}R^{-\frac1m},2^{k}R^{-\frac1m}]\;\;\text{for}\;\; 1\leq k \leq \big[\tfrac{1}{m}\log_{2}R\big].$$
For each $k\geq 1$, we further divide $I_k$  into
\[I_k=\bigcup_{\mu=1}^{2^{\frac{m}{2}(k-1)}} I_{k,\mu}\]
with
$$I_{k,\mu}=
\big[2^{k-1}{R^{-\frac1m}}+(\mu-1)2^{-\frac{m-2}{2}(k-1)}R^{-\frac1m},2^{k-1}R^{-\frac1m}+ \mu 2^{-\frac{m-2}{2}(k-1)}R^{-\frac1m}\big].$$
Thus, we have the following decomposition
\begin{equation}\label{equ:r4regdec}
	[0,1]^{n-1}=\bigcup_{\theta\in \mathcal{F}_{n}(R,m,0,n-1)} \theta,
\end{equation}
where
\begin{align*}
	\mathcal{F}_{n}(R,m,&0,n-1)\\
	:=&\big\{I_{k_1,\mu_1}\times I_{k_2,\mu_2}\times \cdot\cdot\cdot \times I_{k_{n-1},\mu_{n-1}},\;I_0 \times I_{k_2,\mu_2}\times\cdot\cdot\cdot \times I_{k_{n-1},\mu_{n-1}},\cdot\cdot\cdot,\\
	&\quad I_{k_1,\mu_1}\times I_{k_2,\mu_2}\times \cdot\cdot\cdot \times I_0,\;\cdot\cdot\cdot,\;I_0\times I_0\times\cdot\cdot\cdot\times I_0:\;\\
	&\quad 1\leq k_j\leq \big[\tfrac{1}{m}\log_{2}R\big],\;1\leq \mu_j\leq 2^{\frac{m}{2}(k-1)},\;j=1,2,\cdot\cdot\cdot,n-1\big\}.
\end{align*}

\begin{remark} \rm
The $\frac{1}{R}$-neighborhood of the graph of $F^{n-1}_m(0,n-1)$ over each
$\theta \in \mathcal{F}_{n}(R,m,0,n-1)$ can be seen as a rectangular box. For convex hypersurfaces, such a $\theta$ is roughly the largest set near a given point with the same geometric property (Gussian curvature) of the point.
\end{remark}

In terms of the decomposition \eqref{equ:r4regdec}, we refer to \cite{LMZ} for details. It is worth noting that Buschenhenke
\cite{Bus} utilized the analogous decomposition to study the restriction estimates for certain conic hypersurfaces of finite type.

Our main result is the following decoupling inequality
associated with the hypersurface $F^{n-1}_{m}(0,n-1)$ based
on the decomposition \eqref{equ:r4regdec}.

\begin{theorem}
\label{thm:main}
For $2\leq p\leq \frac{2(n+1)}{n-1}$ and each $\varepsilon > 0$, there exists a constant $C(\varepsilon,p)$ such that
\begin{equation}
\label{equ:goalm}
  \|E_{[0,1]^{n-1}}g\|_{L^{p}(B_R)}\leq C(\varepsilon,p)R^{\varepsilon}\Big(\sum_{\theta\in \mathcal{F}_{n}(R,m,0,n-1)}\|E_{\theta}g\|^2_{L^{p}(w_{B_R})}\Big)^{1/2}.
\end{equation}
\end{theorem}

Note that $F^{n-1}_m \vert_{[\frac12,1]^{n-1}}$ has positive definite second fundamental form. We see the sharpness of the exponent in Theorem \ref{thm:main} from the same counterexample of Theorem
\ref{thm:BD}.

%
%
%

The proof of Theorem \ref{thm:main} is based on Theorem \ref{thm:BD} and an induction argument. To make the induction argument complete, we actually deal with a broader class of phase functions $\mathcal{F}^{n-1}$ (the definition can be found in Section 2), which incorporate $F_m^{n-1}(0,m-1)$. To be more precise, we will decompose the cube $[0,1]^{n-1}$ into several parts based on the fact whether the Gaussian curvature of the hypersurface at the associated part is close to zero or not. For the part where the Gaussian curvature is bounded away from zero, it essentially falls into the elliptic setting, we will apply Theorem \ref{thm:BD} directly. For the part where the Gaussian curvature is near zero, we will use the induction argument. Both the dimension and the scale will be inducted.

Theorem \ref{thm:main} can be used to study the problems of the pointwise convergence and the restriction theory in the setting of finite type, one may refer to \cite{LMZ,LiZhaoZhao21}
for more details.

\vskip 0.2in

{\bf Notations:} For nonnegative quantities $X$ and $Y$, we will write $X\lesssim Y$ to denote the estimate $X\leq C Y$ for some $C>0$. If $X\lesssim Y\lesssim X$, we simply write $X\sim Y$.
 Dependence of implicit constants on the power $p$
 or the dimension will be suppressed;
 dependence on additional parameters will be indicated by subscripts.
 For example, $X\lesssim_u Y$ indicates $X\leq CY$ for some $C=C(u)$. We denote $e(t)=e^{2\pi i t}$.
We denote $[x]$ to be the greatest integer not larger than $x$.
We use $B_r(x_0)$ to denote an arbitrary ball centered at $x_0$
 with radius $r$ in $\mathbb{R}^n$ and abbreviate it by $B_r$ in the context.
For any region $\Omega \subset \mathbb{R}^n$,
we denote the characteristic function on $\Omega$ by $\chi_{\Omega}$.
In $\R^n$, we denote $R^{-1/2}\times \cdot\cdot\cdot \times R^{-1/2} \times R^{-1}$-rectangle to be an $R^{-1/2}$-slab in $\mathbb{R}^n$.
Define the Fourier transform on $\mathbb{R}^n$ by
\begin{equation*}
\aligned \widehat{f}(\xi):= \int_{\mathbb{R}^n}e^{- 2\pi ix\cdot \xi}f(x)\,dx,
\endaligned
\end{equation*}
and the inverse Fourier transform by
\begin{equation*}
\aligned \mathcal{F}^{-1}{f}(x):= \int_{\mathbb{R}^n}e^{2\pi ix\cdot \xi}f(\xi)\,d\xi.
\endaligned
\end{equation*}




\section{Proof of the decoupling theorem}\label{sec:2}
First, we make precise the notion of $m$-nondegenerate in this paper.
\begin{definition}
Let $m>2$ be an integer. We say a smooth function $\phi(t)$ is $m$-nondegenerate on the interval $[0,1]$ if there is some constant $C>0$ such that
$$\;0\leq \phi^{(\ell)}(t) \leq C \;\;\text{with}\; \; \phi''(t)\in [C/2,C],\;\text{for}\;0 \leq \ell \leq m, t\in [0,1],$$
and
$$\;\phi^{(\ell)}\equiv 0\;\text{ for }\;\ell\geq m+1.$$
\end{definition}

To prove Theorem \ref{thm:main}, we also need to consider a family of hypersurfaces in $\mathbb{R}^n$ as follows
\[\mathcal{F}^{n-1}:=\Big\{F^{n-1}_{m}(s,n-1-s):\;0 \leq s \leq n-1\Big\},\]
where
\begin{align*}F^{n-1}_{m}(s,n-1-s):= \Big\{(\xi_1,...,&\xi_{n-1},\phi_1(\xi_1)+\cdots
+\phi_s(\xi_s)+\xi_{s+1}^m+\cdots+\xi_{n-1}^m):
\\
&(\xi_1,...,\xi_{n-1}) \in [0,1]^{n-1},\;0 \leq s \leq n-1\Big\}\end{align*}
with each function $\phi_j$ being $m$-nondegenerate.

Let $\mathbf{q}_s$ be the family of all $R^{-1/2}$-cubes in $[0,1]^s$. For each hypersurface $F^{n-1}_{m}(s,n-1-s)$ in the family $\mathcal{F}^{n-1}$, one can write down the corresponding decomposition $\mathcal{F}_{n}(R,m,s,n-1-s)$, where
\begin{align*}
&\quad \mathcal{F}_{n}(R,m,s,n-1-s)\\
&:=\big\{\delta_s\times I_{k_{s+1},\mu_{s+1}}\times I_{k_{s+2},\mu_{s+2}}\times \cdots \times I_{k_{n-1},\mu_{n-1}},\\
&\quad \quad \quad \delta_s\times I_0 \times I_{k_{s+2},\mu_{s+2}}\times\cdots \times I_{k_{n-1},\mu_{n-1}},  \cdots, \\
	&
	\quad \quad \;\; \delta_s\times I_{k_{s+1},\mu_{s+1}}\times I_{k_{s+2},\mu_{s+2}}\times
	 \cdot\cdot\cdot 
	 \times I_0,\;...,\;\delta_s\times I_0\times I_0\times\cdot\cdot\cdot\times I_0:\;\\
	&\quad \quad \;\; \delta_s\in \mathbf{q}_s,\;1\leq k_j\leq \big[\tfrac{1}{m}\log_{2}R\big],\;1\leq \mu_j\leq 2^{\frac{m}{2}(k-1)},\;j=s+1,s+2,...,n-1\big\}.
\end{align*}

We shall establish decoupling inequalities for the class of hypersurfaces $\mathcal{F}^{n-1}$ associated with the corresponding decompositions $\mathcal{F}_{n}(R,m,s,n-1-s)$ with $0\leq s\leq n-1$.

For inequality \eqref{equ:goalm}, we actually prove an stronger version of that.

\begin{proposition}\label{pro:2m}
Suppose that $2\leq p \leq \frac{2(n+1)}{n-1}$. For each $\varepsilon>0$, there exists a constant $C_\varepsilon$ (uniform over all hypersurfaces in $\mathcal F^{n-1}$) such that
\begin{equation}\label{surfacenew}
\| E_{[0,1]^{n-1}}^s g \|_{L^p(B_R)}\leq C_{\varepsilon}R^{\varepsilon}\Big(\sum_{\theta \in \mathcal{F}_n(R,m,s,n-1-s)} \|E_{\theta}g \|^2_{L^p(w_{B_R})}\Big)^{1/2},
\end{equation}
where $E_{[0,1]^{n-1}}^s$ denotes the Fourier extension operator associated with the graph of
\[\sum^s_{i=1}\phi_i(\xi_i)+\sum^{n-1}_{i=s+1}\xi^m_i\]
for $0\leq s\leq n-1$.
\end{proposition}

We shall prove Proposition \ref{pro:2m} by induction on the dimension and the scales. 
The base case is $n=3$, which was proved in \cite{LiZheng21preprint}, 
and the case of small scales  $1\leq r \leq 10$. Our inductive hypothesis is that Proposition \ref{pro:2m} holds for dimension $d$ with $3\leq d \leq n-1$ or radius $1\leq r \leq \frac{R}{2}$.

To prove Proposition \ref{pro:2m}, for a given $s$ and the corresponding hypersurface $F_{m}^{n-1}(s,n-1-s)$, we divide $[0,1]^{n-1}$ into $[0,1]^{n-1}=\bigcup_b \Omega_{(b)}$, where
\[(b):=(b_{s+1},...,b_{n-1}),\;b_j \in \{0, 1\}\]
for $s+1 \leq j \leq n-1$,
\[\Omega_{(b)}:=[0,1]^s\times 
\prod^{n-1}_{j=s+1}J_{b_j},\]
and
\[J_0:= [0, K^{-1/m}],\;J_1:= [K^{-1/m},1].\]
For technical reasons, $K^{-1/m}$, $R^{-1/m}$ and $\big(\tfrac{R}{K}\big)^{-1/m}$ should be dyadic numbers satisfying $1\ll K \ll R^{\varepsilon}$ for any fixed $\varepsilon > 0$. 
Therefore, we choose $K=2^{ms}$  and $R=2^{ml}\; (s,l \in \mathbb{N})$ to be large numbers satisfying $K \approx R^{\varepsilon^{100}}$.

We observe that the hypersurface $F^{n-1}_{m}(s,n-1-s)$
 has positive definite second fundamental form if
  $(\xi_1,...,\xi_{n-1})\in \Omega_{(1,...,1)}$. 
In this region, one can adopt Bourgain-Demeter's decoupling for the perturbed paraboloid.
 While for the regions
$\Omega_{(b)}$ with $b_{s+1}=...=b_{\tilde{s}}=1, b_{\tilde{s}+1}=...=b_{n-1}=0$ for $\;s \leq \tilde{s} \leq n-2$,
we reduce them to the lower dimensional problems. By the Minkowski inequality and Cauchy-Schwarz inequality, we have
\begin{equation}\label{equ:e01intr}
	\| E_{[0,1]^{n-1}}^s g \|_{L^p(B_R)}
	\lesssim 
	\Big(\sum_{b}\| E_{\Omega_{(b)}}^s g \|^2_{L^p(B_R)}\Big)^{1/2}.
\end{equation}
Let $\mathcal{D}_p(R)$ 
denote the optimal constant such that
\begin{equation}\label{equ:defqpr1}
\| E_{[0,1]^{n-1}}^sg \|_{L^p(B_R)}\leq \mathcal{D}_p(R)\Big(\sum_{\theta\in \mathcal{F}_n(R,m,s,n-1-s)}\| E_{\theta}^s g \|^2_{L^{p}(w_{B_R})}\Big)^{1/2}
\end{equation}
for all $1\leq s\leq n-1$ 
and all hypersurfaces in $\mathcal F^{n-1}$.

We are going to treat different regions in different approaches.
First, we estimate the contribution from part $\Omega_{(b)}$ with  $b_{s+1}=...=b_{n-1}=1$.

\subsection{Decoupling for $\Omega_{(1,...,1)}$}
In this subsection, we will establish the decoupling inequality for part $\Omega_{(1,...,1)}$.
We decompose
$$\Omega_{(1,...,1)}=
 \bigcup_{\lambda_{s+1},...,\lambda_{n-1}} \Omega_{\lambda_{s+1},...,\lambda_{n-1}}, \quad
 \Omega_{\lambda_{s+1},...,\lambda_{n-1}}=\bigcup_{\iota_{s+1},...\iota_{n-1}}\tau_{\lambda_{s+1},...,\lambda_{n-1}}^{\iota_{s+1},...\iota_{n-1}},$$
where
\[\Omega_{\lambda_{s+1},...,\lambda_{n-1}}:=[0,1]^s\times\prod^{n-1}_{j=s+1}[\lambda_j, 2\lambda_j],\quad \tau_{\lambda_{s+1},...,\lambda_{n-1}}^{\iota_{s+1},...\iota_{n-1}}:=\alpha_s\times \prod^{n-1}_{j=s+1}J_{\lambda_j,\iota_j}.\]
Here $\alpha_s\subset [0,1]^s$ are $K^{-1/2}$-cubes and
\[J_{\lambda_j,\iota_j}:=[\lambda_j+(\iota_j-1)\lambda_j^{-\frac{m-2}{2}}K^{-1/2}, \lambda_j+\iota_j\lambda_j^{-\frac{m-2}{2}}K^{-1/2}]\]
for $1\leq \iota_j \leq \lambda_j^{1+\frac{m-2}{2}}K^{1/2}$ and $\lambda_j\in [K^{-1/m}, \frac{1}{2}]$ is a dyadic number.
In the following, we will abbreviate $\tau_{\lambda_{s+1},...,\lambda_{n-1}}^{\iota_{s+1},...\iota_{n-1}} $ as $\tau$ for convenience.

Combining this decomposition with the Minkowski inequality and Cauchy-Schwarz inequality, we get a trivial decoupling at the scale $K$ for
$2\leq p \leq \frac{2(n+1)}{n-1}$
\begin{equation}\label{equ:trivial}
\| E_{\Omega_{(1,...,1)}}^sg \|_{L^p(B_K)}\lesssim K^{\frac{n-1}{4}}\Big(\sum_{\tau \subset \Omega_{(1,...,1)}} \| E_{\tau}^s g \|^2_{L^p(B_K)}\Big)^{1/2}.
\end{equation}
Summing over all the balls $B_K \subset B_R$, we obtain
\begin{equation}\label{equ:trivialsumming}
\| E_{\Omega_{(1,...,1)}}^sg \|_{L^p(B_R)}\lesssim K^{\frac{n-1}{4}}\Big(\sum_{\tau \subset \Omega_{(1,...,1)}} \| E_{\tau}^s g \|^2_{L^p(B_R)}\Big)^{1/2}.
\end{equation}
We claim that, for any given $\tau\subset\Omega_{(1,...,1)}$ of size $K^{-1/2}\times...\times K^{-1/2}\times \lambda_{s+1}^{-\frac{m-2}{2}}K^{-1/2}\times...\times \lambda^{-\frac{m-2}{2}}_{n-1}K^{-1/2}$, we have
\begin{lemma}\label{lem:omega0tau}
For $2\leq p \leq \frac{2(n+1)}{n-1}$ and each $\varepsilon > 0$, there exists a positive constant $C_\varepsilon$ such that
\begin{equation}\label{equ:lamsigma}
  \|E_{\tau}^s g\|_{L^p(B_R)}
  \leq C_{\varepsilon}\big(\frac{R}{K}\big)^{\varepsilon}\Big(\sum_{\theta \subset \tau} \| E_{\theta}^s g \|^2_{L^p(w_{B_R})}\Big)^{1/2},
\end{equation}
where $\theta \in \mathcal{F}_n(R,m,0,n-1)$.
\end{lemma}

With Lemma \ref{lem:omega0tau} in hand, plugging \eqref{equ:lamsigma} into \eqref{equ:trivialsumming}, we get the decoupling inequality for $\Omega_{(1,...,1)}$
\begin{equation}\label{equ:omegaoo}
\| E_{\Omega_{(1,...,1)}}^s g \|_{L^p(B_R)}\leq C_{\varepsilon}K^{O(1)}R^{\varepsilon}\Big(\sum_{\theta \subset \Omega_{(1,...,1)}} \Vert E_{\theta}^s g \Vert^2_{L^p(w_{B_R})}\Big)^{1/2}.
\end{equation}

\begin{proof}[{\bf Proof of Lemma \ref{lem:omega0tau}:}]
By change of variables, we may assume that $B_R$ is centered at the origin and
\[\tau = [0, K^{-1/2}]^s\times \prod^{n-1}_{j=s+1}[\lambda_j, \lambda_j+\lambda_j^{-\frac{m-2}{2}}K^{-1/2}].\]
By a change of variables, we see that
\[\|E_{\tau}^s g\|^p_{L^p(B_R)}=(\lambda_{s+1}^{\frac{m-2}{2}}\cdot\cdot\cdot \lambda_{n-1}^{\frac{m-2}{2}})^{-p+1}K^{\frac{n+1}{2}-\frac{n-1}{2}p} \|E^{S}_{[0,1]^{n-1}}\tilde{g}\|^p_{L^p(\mathcal{L}_0(B_R))},\]
where $E_{[0,1]^{n-1}}^S$ is the same operator defined as in Theorem \ref{thm:BD}
 and $\mathcal{L}_0(B_R)$
denotes the image of $B_R$ with size roughly as follows
\[K^{-1/2}R\times\cdots\times K^{-1/2}R \times \lambda_{s+1}^{-\frac{m-2}{2}}K^{-1/2}R\times\cdots \times \lambda^{-\frac{m-2}{2}}_{n-1}K^{-1/2}R\times K^{-1}R.\]
In fact, $S$ is the graph of
\[
\begin{split}
\tilde
\psi( \tilde{\xi}_1,...,\tilde{\xi}_{n-1})
& =
\sum_{j=1}^{s} \Big(K\phi_j(K^{-\frac12} \tilde{\xi}_j )-K\phi_j(0)-K^{\frac{1}{2}}\phi'_j(0)\xi_j \Big)\\
&
\quad +
 \sum_{j=s+1}^{n-1}  \Big(  (\lambda_j^{-\frac{m-2}{2}} K^{-\frac12}\tilde{\xi}_{j} + \lambda_j)^m K - \lambda_j^{m} K - (m-1) \lambda_j^{\frac{m}{2}} K^{\frac12}\tilde{\xi}_{j}  \Big)
\end{split}
\]
on $[0,1]^{n-1}$. It is easy to see that $S$ has positive definite second fundamental form.

We divide $\mathcal{L}_0(B_R)$  into a finitely overlapping union of balls as follows
\[\mathcal{L}_0(B_R)=\bigcup B_{R/K}.\]
For a given $\theta \subset \tau$ such as
\[\theta = [0, R^{-1/2}]^s\times[\lambda_{s+1}, \lambda_{s+1}+\lambda^{-\frac{m-2}{2}}_{s+1}R^{-1/2}]\times\cdots \times[\lambda_{n-1}, \lambda_{n-1}+\lambda^{-\frac{m-2}{2}}_{n-1}R^{-1/2}],\]
under a change of variables we deduce that the image of $\theta$ is
\[\tilde{\theta}=[0, K^{1/2}R^{-1/2}]^{n-1}.\]
We use Bourgain-Demeter's decoupling inequality \eqref{equ:bd}
 on each $B_{R/K}$ to obtain
\begin{equation}
\label{omega0tauparp}
\big\|E_{[0,1]^2}^S
 \tilde{g}\big\|_{L^p(B_{R/K})}\leq C_{\varepsilon}\big(\frac{R}{K}\big)^{\varepsilon}\Big(\sum_{\tilde{\theta}:K^{1/2}R^{-1/2}-cube}\| E_{\tilde \theta}^S  \tilde{g}\|^2_{L^p(w_{B_{R/K}})}\Big)^{1/2}.
\end{equation}

Summing over all the balls $B_{R/K}\subset \mathcal{L}_0(B_R)$ on both sides of \eqref{omega0tauparp} and using Minkowski's inequality, we have
\[\|E^{S}_{[0,1]^{n-1}}\tilde{g}\|_{L^p(\mathcal{L}_0(B_R))}
\leq C_{\varepsilon}
\big(\frac{R}{K}\big)^{\varepsilon}\Big(\sum_{\tilde{\theta}:K^{1/2}R^{-1/2}-cube}\|
E^{S}_{\tilde{\theta}}\tilde{g}
\|^2_{L^p(\tilde{w}_{\mathcal{L}_0(B_R)})}
\Big)^{1/2},\]
where 
$\tilde{w}_{\mathcal{L}_0(B_R)} :=\sum\limits_{B_{R/K} \subset \mathcal{L}_0 (B_R)} w_{B_{R/K}}.$
Taking the inverse change of variables, it follows
\[
\|E_{\tau}^sg\|_{L^p(B_R)}\leq C_{\varepsilon}\big(\frac{R}{K}\big)^{\varepsilon}\Big(\sum_{\theta \subset \tau}\| E_{\theta}^sg\|^2_{L^p(w_{B_R})}\Big)^{1/2}.
\]
Therefore, we complete the proof of Lemma \ref{lem:omega0tau}.
\end{proof}

\subsection{Decoupling for $\Omega_{(b)}\;(b_{s+1}=\cdots=b_{\tilde{s}}=1,\;b_{\tilde{s}+1}=\cdots=b_{n-1}=0,\;s \leq \tilde{s}\leq n-2)$}

We decompose
$$\Omega_{(b)}=
 \bigcup_{\mu_{s+1},...,\mu_{n-1}} \Omega_{\mu_{s+1},...,\mu_{n-1}}, \quad
 \Omega_{\mu_{s+1},...,\mu_{n-1}}=\bigcup_{\nu_{s+1},...,\nu_{n-1}}\tau_{\mu_{s+1},...,\mu_{n-1}}^{\nu_{s+1},...,\nu_{n-1}},$$
where
\[\Omega_{\mu_{s+1},...,\mu_{n-1}}:=[0,1]^s\times\prod^{\tilde{s}}_{j=s+1}[\mu_j, 2\mu_j]\times [0,K^{-1/m}]^{n-1-\tilde{s}}\]
and \[ \tau_{\mu_{s+1},...,\mu_{n-1}}^{\nu_{s+1},...,\nu_{n-1}}:=\beta_s\times \prod^{\tilde{s}}_{j=s+1}J_{\mu_j,\nu_j}\times [0,K^{-1/m}]^{n-1-\tilde{s}}.\]
Here $\beta_s\subset [0,1]^s$ are $K^{-1/2}$-cubes and
\[J_{\mu_j,\nu_j}:=[\mu_j+(\nu_j-1)\mu_j^{-\frac{m-2}{2}}K^{-1/2}, \mu_j+\nu_j\mu_j^{-\frac{m-2}{2}}K^{-1/2}]\]
for $1\leq \nu_j \leq \mu_j^{1+\frac{m-2}{2}}K^{1/2}$ and $\mu_j\in [K^{-1/m}, \frac{1}{2}]$ is a dyadic number.
In the following, we will abbreviate $\tau_{\mu_{s+1},...,\mu_{n-1}}^{\nu_{s+1},...,\nu_{n-1}} $ as $\tau$ for convenience.

First, by reduction of dimension arguments,
we are able to prove the following result.
\begin{lemma}\label{lem:3.3}
For $2\leq p\leq \frac{2n}{n-2}$ and any $\varepsilon>0$, there exists a constant $C_\varepsilon$ such that
\begin{equation}\label{omegall}
\| E_{\Omega_{(b)}}^{s}g \|_{L^p(B_R)}\leq C_{\varepsilon}K^{\varepsilon}\Big(\sum_{\tau \subset \Omega_{(b)}} \| E_{\tau}^{s}g \|^2_{L^p(w_{B_R})}\Big)^{1/2}.
\end{equation}
\end{lemma}

To prove Lemma \ref{lem:3.3}, we employ the induction hypothesis on the dimension. Freeze the $x_{n-1}$ variable and fix a bump function $\varphi \in C^{\infty}_c(\mathbb{R}^n)$ with ${\rm supp}\; \varphi \subset B(0,1)$
 and $\vert \check{\varphi}(x)\vert \geq 1$ for all $x \in B_1(0)$.
Define $F:= \check{\varphi_{K^{-1}}}E_{\Omega_{(b)}}g$, where $\varphi_{K^{-1}}(\xi):= K^n \varphi(K\xi),\; \xi \in \mathbb{R}^n$.
We denote $F(\cdot,x_{n-1},\cdot)$ by $G$.
It is easy to see that ${\rm supp}\;\hat{G}$ is contained in the projection of ${\rm supp}\;\hat{F}$ on the hyperplane $\xi_{n-1}=0$, that is,
 in the $K^{-1}$-neighborhood of $F^{n-2}_{m}(s,n-2-s)$.
 By the induction hypothesis on the dimension, we have
\[
\|G\|_{L^p(\mathbb{R}^{n-1})}\leq C_{\varepsilon}K^{\varepsilon} \Big(\sum_{\tau \subset \Omega_{(b)}}\Vert G_{\tau}\Vert^2_{L^p(\mathbb{R}^{n-1})}\Big)^{1/2},
\]
i.e.,
\[
\begin{split}
\quad  \|F(\cdot,x_{n-1},\cdot)\|_{L^p(\mathbb{R}^{n-1})}
 \leq
C_{\varepsilon} K^{\varepsilon} \Big(\sum_{\tau \subset \Omega_{(b)}}\Vert F_{\tau}(\cdot,x_{n-1},\cdot)\Vert^2_{L^p(\mathbb{R}^{n-1})}\Big)^{1/2},
\end{split}
\]
where $F_{\tau}(x):=\check{\varphi_{K^{-1}}}E_{\tau}^s g$.
Integrating on both sides of the above inequality with respect to $x_{n-1}$-variable from
 $-\infty$ to $\infty$, 
 we derive
\[\|F\|_{L^p(\mathbb{R}^n)}\leq C_{\varepsilon}K^{\varepsilon} \Big(\sum_{\tau \subset \Omega_{(b)}}\Vert F_{\tau}\Vert^2_{L^p(\mathbb{R}^n)}\Big)^{1/2}.\]
Thus, we have
\begin{align*}
	\| E_{\Omega_{(b)}}^{s}g \|_{L^p(B_K)}\lesssim
	 C_{\varepsilon}K^{\varepsilon}
	 \Big(\sum_{ \tau \subset \Omega_{(b)}} \| E_{\tau}^{s}g \|^2_{L^p(w_{B_K})}\Big)^{1/2}.
\end{align*}
Summing over all the balls $B_K\subset B_R$, we get the inequality \eqref{omegall} as required.

Next, we estimate the $\| E_{\tau}^{s}g \|_{L^p(B_R)}$ on the right-hand side of \eqref{omegall}.  For each $\tau \subset \Omega_{b}$, we claim that
\begin{equation}\label{omegas1}
	\| E_{\tau}^{s}g \|_{L^p(B_R)}
	\leq
	C_{\varepsilon}\big(\tfrac{R}{K}\big)^{\varepsilon}
	\Big(
	\sum_{\theta \subset \tau}
	\| E_{\theta}^{s}g \|^2_{L^p(w_{B_R})}
	\Big)^{1/2},\quad 2\leq p \leq \frac{2(n+1)}{n-1},
\end{equation}
where $\theta\in\mathcal{F}_n(R,m,s,n-1-s)$. 
We take the change of variables associated with $\tau$ such that $\vert E_{\tau}^s g(x)\vert$ becomes $\vert E_{[0,1]^{n-1}}^{\tilde{s}}\tilde{g}(\tilde{x})\vert$ 
for some 
$s\leq \tilde{s} \leq n-2$
 and some function $\tilde{g}(\tilde{x})$ with $\tilde{x}\in \mathcal{L}_{\lambda}(B_R)$, where $\mathcal{L}_{\lambda}(B_R)$ denotes the image of $B_R$.
  We rewrite $\mathcal{L}_{\lambda}(B_R)$ into a finitely overlapping union of balls as follows
\[\mathcal{L}_{\lambda}(B_R)=\bigcup B_{R/K}.\]

We will use an induction argument to complete the proof. To do this, we need verify that for $\theta \subset \tau$, its image $\tilde{\theta}$ belongs to the decomposition $\mathcal{F}_n(\frac{R}{K},m;\tilde{s},n-1-\tilde{s})$ under the change of variables associated with $\tau$. Without loss of generality, we may assume that
\[
\tau=[0,K^{-1/2}]^s\times \prod^{\tilde{s}}_{l=s+1}[\mu_l,\mu_l+\mu_l^{-\frac{m-2}{2}}K^{-1/2}]\times [0,K^{-1/m}]^{n-1-\tilde{s}}
\]
and
\[
\theta=[0,R^{-1/2}]^s
\times \prod^{s'}_{k=s+1}[\mu_k, \mu_k+\mu_k^{-\frac{m-2}{2}}R^{-1/2}]\times[0, R^{-1/m}]^{n-1-s'},
\]
where $\tilde{s}\leq s' \leq n-2$. Under the change of variables associated with $\tau$, we see that
\begin{align*}
\tilde{\theta}& =
[0,(\frac{R}{K})^{-1/2}]^{s}
\times \prod_{k=s+1}^{\tilde s} [\mu_k^\frac{m}2K^\frac12,\mu_k^\frac{m}2K^\frac12+(\frac{R}K)^\frac12]\\
& \quad\quad\quad\quad\quad\quad\quad \times \prod^{s'}_{k=\tilde{s}+1}[\tilde{\mu}_k, \tilde{\mu}_k+\tilde{\mu}_k^{-\frac{m-2}{2}}(\frac{R}{K})^{-1/2}]\times[0,(\frac{R}{K})^{-1/m}]^{n-1-s'},
\end{align*}
where $\tilde{\mu}_k:= K^{\frac{1}{m}}\mu_k$ is also a dyadic number. 
Thus, $ \tilde \theta \in \mathcal{F}_n\big(\tfrac{R}{K},m,\tilde s,n-1- \tilde s\big).$

By our induction hypothesis, we have
\[\|E_{[0,1]^{n-1}}^{\tilde{s}} \tilde{g}\|_{L^p(B_{R/K})}\leq C_{\varepsilon}\big(\tfrac{R}{K}\big)^{\varepsilon}\Big(\sum_{\tilde\theta\in \mathcal{F}_n(\frac{R}{K},m,\tilde{s},n-1-\tilde{s})}\|E_{\tilde{\theta}}^{\tilde{s}} \tilde{g}\|^2_{L^p(w_{B_{R/K}})}\Big)^{1/2},\]
for $2\leq p \leq \frac{2(n+1)}{n-1}$.
Summing over all the balls $B_{R/K}\subset \mathcal{L}_{\lambda}(B_R)$ and using Minkowski's inequality, we have
\[
\|E_{[0,1]^{n-1}}^{\tilde{s}} \tilde{g}\|_{L^p(\mathcal{L}_{\lambda}(B_R))}
\leq
 C_{\varepsilon}
 \big(\tfrac{R}{K}\big)^{\varepsilon}\Big(\sum_{\tilde\theta\in \mathcal{F}_n(\frac{R}{K},m,\tilde s,n-1-\tilde s)}\| E_{\tilde{\theta}}^{\tilde{s}}\tilde{g}\|^2_{L^p(\tilde{w}_{\mathcal{L}_{\lambda}(B_R)})}\Big)^{1/2},
\]
where 
$\tilde{w}_{\mathcal{L}_{\lambda}(B_R)} :=\sum\limits_{B_{R/K} \subset \mathcal{L}_{\lambda}(B_R)}w_{B_{R/K}}.$
Taking the inverse change of variables, we deduce
\[\|E_{\tau}^s g\|_{L^p(B_R)}\leq C_{\varepsilon}\big(\tfrac{R}{K}\big)^{\varepsilon}\Big(\sum_{\theta \subset \tau}\| E_{\theta}^s g\|^2_{L^p(w_{B_R})}\Big)^{1/2},\]
and thus Claim \eqref{omegas1} is verified.

Plugging \eqref{omegas1} into \eqref{omegall}, one has
\begin{align}\nonumber
\| E_{\Omega_{(b)}}^{s}g \|_{L^p(B_R)}\leq& C_{\varepsilon}K^{\varepsilon}\big(\tfrac{R}{K}\big)^{\varepsilon}\Big(\sum_{\tau \subset \Omega_{(b)}}\sum_{\theta \subset \tau} \|E_{\theta}^{s}g \|^2_{L^p(w_{B_R})}\Big)^{1/2}\\\label{equ:controme1est}
=& C_{\varepsilon}R^{\varepsilon}\Big(\sum_{\theta \subset \Omega_{(b)}} \| E_{\theta}^s g \|^2_{L^p(w_{B_R})}\Big)^{1/2},
\end{align}
where $\theta \in \mathcal{F}_n(R,m,s,n-1-s)$.
Thus, we obtain the decoupling inequality for $\Omega_{(b)}$-part. 

In summary, we get
\[\mathcal{D}_p(R)\leq C(\varepsilon')K^{O(1)}R^{\varepsilon'}+2^{n-1}C_{\varepsilon''}K^{\varepsilon''}C_{\varepsilon}(\frac{R}{K})^{\varepsilon}.\]
For any given $\varepsilon>0$, we choose $\varepsilon'=\varepsilon''=\varepsilon^2$. Recall that $K\approx R^{\varepsilon^{100}}$. We conclude that there exists a constant $R_0$ depending only on $\varepsilon$ such that
\[\mathcal{D}_p(R)\leq C_{\varepsilon}R^{\varepsilon}\]
whenever $R>R_0$.
This completes the proof of Proposition \ref{pro:2m}.

\begin{remark}
It is also an interesting problem to investigate whether the argument in this paper can be applied to study the model in \cite{LiYang21, Kemp, LiYang21preprint}. It is worth noting that the sharp decoupling inequalities for general smooth hypersurfaces with vanishing curvature in
$\mathbb{R}^n$ for $n \geq 4$ are still open.
\end{remark}





\bibliographystyle{amsplain}


\end{document}